\documentclass[11pt]{article}
\usepackage{amssymb,amsmath}

\usepackage{graphics}
\newtheorem{theorem}{Theorem}

\def\cadre{$$\vcenter\bgroup\advance\hsize by -4em\noindent%
\ignorespaces\it\refstepcounter{equation}}\makeatletter
\def\endcadre{\rm\egroup\leqno(\theequation)$$\global\@ignoretrue}
\makeatother

\begin{document}

\title{Graphs with no induced wheel or antiwheel}

\author{Fr\'ed\'eric Maffray\thanks{Partially supported by ANR project
STINT under reference ANR-13-BS02-0007.} \\
{\small CNRS, Laboratoire G-SCOP, Universit\'e de Grenoble, 
France}  \\ {\small frederic.maffray@grenoble-inp.fr}}

\date{\today}

% \institute{\today}

% \titlerunning{ } \maketitle
%

\maketitle
\begin{abstract}
A wheel is a graph that consists of a chordless cycle of length at
least~$4$ plus a vertex with at least three neighbors on the cycle.
It was shown recently that detecting induced wheels is an NP-complete
problem.  In contrast, it is shown here that graphs that contain no
wheel and no antiwheel have a very simple structure and consequently
can be recognized in polynomial time.
\end{abstract}

% \section{Introduction}

Four families of graphs have repeatedly played important roles in
structural graph theory recently.  They are called \emph{Truemper
configurations} as they were first used by Truemper in several
theorems \cite{Tru}.  These configurations are called \emph{pyramids},
\emph{prisms}, \emph{thetas} and \emph{wheels}.  We will not recall
all the definitions, as we do not need all of them here; see
Vuskovi\'c \cite{Vus} for a very extensive survey on Truemper
configurations.  It is interesting to know the complexity of deciding
whether a graph contains a Truemper configuration of a certain type.
The problem is polynomial for pyramids \cite{CCLSV}; indeed it is one
of the main steps in the polynomial-time recognition algorithm for
perfect graphs presented in \cite{CCLSV}.  On the other hand, the
problem is NP-complete for thetas \cite{CK} and prisms \cite{MafTro}.
Here we will deal only with the fourth Truemper configuration, the
wheel.  A \emph{wheel} is a graph that consists of a cycle of length
at least~$4$ plus a vertex that has at least three neighbors on the
cycle.  Diot, Tavenas and Trotignon \cite{DTT} proved that it is also
NP-complete for wheels, and they mention the question of
characterizing the graphs that contain no wheel and no antiwheel but
leaves it open.  This question is solved here with a complete
description of the structure of these graphs, from which it follows
easily that they can be recognized in polynomial (indeed linear) time.

We use the standard graph-theoretic terminology.  We let $K_n$, $P_n$
and $C_n$ respectively denote the complete graph, path and cycle on
$n$ vertices, and $nF$ denote the graph with $n$ components, all
isomorphic to $F$.  Given a graph family $\cal{F}$, a graph $G$ is
\emph{$\cal{F}$-free} if no induced subgraph of $G$ is isomorphic to
any member of $\cal{F}$; when $\cal{F}$ has only one element $F$ we
say that $G$ is \emph{$F$-free}.

We recall the following simple characterization of $P_5$-free
bipartite graphs.  
\begin{theorem}[See {\cite{HPS}},
{\cite[Section~2.4]{MahPel}}]\label{thm:bp5}
Let $H$ be a connected bipartite graph, where $V(H)$ is partitioned
into stable sets $X$ and $Y$.  The following conditions are
equivalent:
\vspace{-.3cm}
\begin{itemize}
\item 
$H$ is $P_5$-free; 
\item 
$H$ is $2K_2$-free; 
\item 
The neighborhoods of any two vertices in $X$ are comparable by
inclusion (equivalently, the same holds in $Y$); 
\item 
There is an integer $h>0$ such that $X$ can be partitioned into
non-empty sets $X_1, \ldots, X_h$ and $Y$ can be partitioned into
non-empty sets $Y_1, \ldots, Y_h$ such that for all $i,j\in\{1,
\ldots, h\}$ a vertex in $X_i$ is adjacent to a vertex in $Y_j$ if and
only if $i+j\le h+1$.
\end{itemize}
\end{theorem}
It follows from Theorem~\ref{thm:bp5} that when $H$ is a $P_5$-free
connected bipartite graph, with the same notation as in the theorem,
then $X$ contains a vertex that is complete to $Y$, and $Y$ contains a
vertex that is complete to $X$.

Recall that a graph is \emph{split} \cite{FH} if its vertex-set can be
partitioned into a stable set and a clique.  F\"oldes and Hammer
\cite{FH} gave the following characterization of split graphs.
\begin{theorem}[\cite{FH}]\label{thm:FH}
A graph is split if and only if it is $\{2K_2, C_4, C_5\}$-free.
\end{theorem}

We define three classes of graphs $\cal A$, $\cal B$ and $\cal C$ as
follows.

\noindent{\bf Class $\cal A$:} A graph $G$ is $\cal A$ if $V(G)$ can
be partitioned into two sets $\{a,b,c,d,e\}$ and $X$ such that:
\vspace{-.3cm}
\begin{itemize}
\setlength{\itemsep}{0pt}
\item
$\{a,b,c,d\}$ induces a hole with edges $ab, bc, cd, da$; 
\item
$X$ is non-empty, induces a clique and is complete to $\{c,d\}$ and
anticomplete to $\{a,b\}$;  
\item 
$e$ is complete to $X$, anticomplete to $\{a,b\}$, and has at most one
neighbor in $\{c,d\}$.
\end{itemize}

\noindent{\bf Class $\cal B$:} A graph $G$ is in $\cal B$ if $V(G)$
can be partitioned into four stable sets $X, Y, Z, W$, with two
special vertices $x\in X$ and $y\in Y$, such that: %
\vspace{-.3cm}
\begin{itemize}
\setlength{\itemsep}{0pt}
\item
$|X|\ge 2$, $|Y|\ge 2$, and $X\cup Y$ induces a connected $P_5$-free
bipartite graph;
\item
$W$ is anticomplete to $X\cup Y\cup Z$ (so all vertices of $W$ are
isolated in $G$);
\item
$x$ is complete to $Y$ and $y$ is complete to $X$;
\item
$Z$ is complete to $\{x,y\}$ and anticomplete to $(X\cup
Y)\setminus\{x,y\}$.
\end{itemize}

\noindent{\bf Class $\cal C$:} A graph $G$ is in $\cal C$ if $V(G)$
can be partitioned in two cliques $X$ and $Y$ of size at least $2$
such that the edges between $X$ and $Y$ form a matching of size~$2$.

\begin{theorem}\label{thm:main}
The following three properties are equivalent:
\vspace{-.3cm}
\begin{itemize}
\setlength{\itemsep}{0pt}
\item 
$G$ is (wheel, antiwheel)-free.
\item
$G$ contains no wheel or antiwheel on at most seven vertices.
\item
$G$ or $\overline{G}$ is either a $5$-hole, a $6$-hole, a split graph,
or a member of ${\cal A}\cup {\cal B}\cup {\cal C}$.
\end{itemize}
\end{theorem}
\noindent{\it Proof.} Let $F_1$ (resp.~$F_2$) be the wheel that
consists of a $4$-hole plus a vertex adjacent to three (resp.~four)
vertices of the hole.  

Clearly, the first condition of the theorem implies the second.
Suppose that $G$ satisfies the third condition.  If $G$ or
$\overline{G}$ is a $5$-hole or a $6$-hole, then clearly it does not
contains a wheel.  If $G$ is a split graph, it contains no hole and
consequently no wheel.  Suppose that $G\in {\cal A}\cup {\cal B}\cup
{\cal C}$.  If $G\in {\cal A}\cup {\cal C}$, it contains only one hole
$H$, of length~$4$.  If $G\in {\cal B}$ it may contain many holes, but
they all have four vertices, more precisely two vertices from $X$ and
two from $Y$.  In all cases, it is easy to see that whenever $H$ is a
hole in $G$, every vertex of $G\setminus H$ has at most two neighbors
in $H$.  So no hole of $G$ extends to a wheel.

Now let us prove that the second condition implies the third.  Let $G$
be a graph that contains no wheel or antiwheel on at most seven
vertices.

First suppose that $G$ contains a $5$-hole $C$.  Note that $V(C)$ also
induces a $5$-hole in $\overline{G}$.  If there is any vertex $x$ in
$V(G)\setminus V(C)$, then $x$ has either at least three neighbors in
$C$ or three non-neighbors in $C$, and so $V(C)\cup\{x\}$ induces a
wheel in $G$ or in $\overline{G}$.  Thus no such $x$ exists, and $G$
is a $5$-hole.

Now suppose that $G$ contains a $6$-hole $C$, with vertices $c_1,
\ldots, c_6$ and edges $c_ic_{i+1}$, with subscripts modulo $6$.  Pick
any $x$ in $V(G)\setminus V(C)$.  Vertex $x$ has at most two neighbors
in $C$, for otherwise $(C,x)$ is a wheel in $G$.  It follows that, up
to symmetry, $N(x)\cap V(C)$ is equal either to $\{c_1\}$, $\{c_1,
c_2\}$, or $\{c_1, c_5\}$, and in that case $\{x, c_1, c_3, c_4,
c_6\}$ induces an $\overline{F}_1$, or to $\{c_1, c_4\}$, and in that
case $\{x, c_2, c_3, c_5, c_6\}$ induces an $\overline{F}_2$.  Thus no
such $x$ exists, and $G$ is a $6$-hole.

If $G$ contains a $6$-antihole, then the same argument as in the
preceding paragraph, applied to $\overline{G}$, implies that $G$ is a
$6$-antihole.

Now assume that $G$ contains no $5$-hole, no $6$-hole and no
$6$-antihole.  We may also assume that $G$ is not a split graph, for
otherwise the theorem holds.  It follows from Theorem~\ref{thm:FH}
that $G$ contains either a $2K_2$, a $C_4$ or a $C_5$.  Since $G$
contains no $C_5$, and up to self-complementation, we may therefore
assume that $G$ contains a $2K_2$.  Let $A,B$ be two disjoint subsets
of $V(G)$ such that both $A$ and $B$ are cliques of size at least $2$
and $A$ is anticomplete to $B$.  Graph $G$ admits such a pair since we
can let $A$ and $B$ be the two cliques of size $2$ of a $2K_2$.
Choose $A$ and $B$ such that $|A\cup B|$ is maximized.  Let
$R=V(G)\setminus (A\cup B)$.  We observe that:
\begin{cadre}\label{R1}
For every vertex $x$ in $R$, either: \\
$\bullet$ $x$ is complete to $A$ and has a neighbor in $B$, or
vice-versa, or \\
$\bullet$ $x$ has exactly one non-neighbor in $A$ and one non-neighbor
in~$B$.
\end{cadre}
Indeed, if $x$ has at most one non-neighbor in $A$ and at most one
non-neighbor in $B$, then (\ref{R1}) holds.  So suppose, up to
symmetry, that $x$ has two non-neighbors $a,a'$ in $A$.  If $x$ has a
non-neighbor $b$ in $B$, then, picking any $b'\in B\setminus b$, we
see that $\{x, a, a', b, b'\}$ induces an $\overline{F}_1$ or
$\overline{F}_2$ (depending on the pair $x,b'$), a contradiction.  So
$x$ is complete to $B$.  If $x$ has no neighbor in $A$, then the pair
$A,B\cup\{x\}$ contradicts the choice of $A,B$.  So $x$ has a neighbor
in $A$, and the first item in (\ref{R1}) holds.  This proves
(\ref{R1}).

Let $A=\{a_1, \ldots, a_p\}$, with $p\ge 2$, and let $B=\{b_1, \ldots,
b_q\}$, with $q\ge 2$.  Define the following subsets of $R$: \\
$\bullet$ $R_0 = \{x\in R \mid x$ is complete to $A$ or to $B\}$.  \\
$\bullet$ $R_{i,j} = \{x \in R \mid x$ is complete to $(A\cup B)\setminus
\{a_i, b_j\}$ and anticomplete to $\{a_i, b_j\}\}$, for each $(i,j)\in
\{1, \ldots, p\}\times \{1, \ldots, q\}$.  \\
Clearly these sets are pairwise disjoint, and (\ref{R1}) means that
$R=R_0 \cup \bigcup_{i,j} R_{i,j}$.  

\medskip

Say that two vertices $x$ and $y$ of $R$ are \emph{$A$-comparable} if
one of the two sets $N_A(x)$ and $N_A(y)$ contains the other; in the
opposite case, say that $x$ and $y$ are \emph{$A$-incomparable}.
Define the same with respect to $B$.

Suppose that there are two $A$-incomparable vertices $x$ and $y$ in
$R$.  Up to relabeling, $a_1$ is adjacent to $x$ and not to $y$ and
$a_2$ is adjacent to $y$ and not to $x$.  Since each of $x$ and $y$
has a neighbor in $B$, there is a chordless path $P$ whose endvertices
are $x$ and $y$ and whose interior vertices are in $B$; and since $B$
is a clique, the length $\ell$ of $P$ is equal to $2$ or $3$.  We may
assume that if $\ell=2$ then $P=x$-$b_1$-$y$ while if $\ell=3$ then
$P=x$-$b_1$-$b_2$-$y$.  Vertices $x$ and $y$ are adjacent, for
otherwise $V(P)\cup\{a_1, a_2\}$ induces a $5$-hole or a $6$-hole.
Note that if $p\ge 3$, then $x$ has no neighbor $a$ in
$A\setminus\{a_1, a_2\}$, for otherwise $\{a_1, a_2, x, y, a\}$
induces an $F_1$ or $F_2$; and the same holds for $y$.  So if $p\ge
3$, $x$ and $y$ are anticomplete to $A\setminus\{a_1, a_2\}$ and, by
(\ref{R1}), they are complete to $B$.

Let $z$ be any vertex in $R\setminus \{x, y\}$.  Suppose that $z$ is
complete to $\{a_1, a_2\}$.  Then $z$ is anticomplete to $\{x,y\}$,
for otherwise $\{x, y, z, a_1, a_2\}$ induces an $F_1$ or $F_2$.  Then
$z$ is not adjacent to $b_1$, for otherwise either $\{x, y, z, b_1,
a_1, a_2\}$ induces a $6$-antihole (if $\ell=2$) or $\{x, y, z, b_1,
a_2\}$ induces a $5$-hole (if $\ell=3$).  Let $b$ be a neighbor of $z$
in $B$; so $b\neq b_1$.  Then $x$ is adjacent to $b$, for otherwise
$\{x, z, b, b_1, a_1\}$ induces a $5$-hole, and $y$ is adjacent to
$b$, for otherwise $\{x, y, z, b, a_2\}$ induces a $5$-hole; but then
$\{x, y, z, b, a_1, a_2\}$ induces a $6$-antihole.  It follows that no
vertex of $R$ is complete to $\{a_1, a_2\}$.  By the same argument, if
$\ell=3$ then no vertex of $R$ is complete to $\{b_1, b_2\}$, and
consequently $R_0=\emptyset$.

Suppose that $\ell=3$.  The preceding arguments and (\ref{R1}) imply
that $R=R_{1,1}\cup R_{1,2}\cup R_{2,1}\cup R_{2,2}$.  Note that $x\in
R_{2,2}$ and $y\in R_{1,1}$.  If $p\ge 3$, then $\{x, y, a_1, a_2,
a_3\}$ induces an $F_2$.  So $p=2$, and similarly $q=2$.  If there is
any vertex $u$ in $R_{1,2}$, then $u$ is adjacent to $x$, for
otherwise $\{u, x, a_1, a_2, b_1\}$ induces a $5$-hole, and similarly
$u$ is adjacent to $y$; but then $\{u, x, y, a_1, a_2\}$ induces an
$F_1$.  So $R_{1,2}=\emptyset$, and similarly $R_{2,1}=\emptyset$.
Therefore $V(G)=\{a_1, a_2, b_1, b_2\}\cup R_{1,1}\cup R_{2,2}$.  If
some vertex $u$ in $R_{1,1}$ is not adjacent to some vertex $v$ in
$R_{2,2}$, then $\{u, v, a_1, a_2, b_1, b_2\}$ induces a $6$-hole.  So
$R_{1,1}$ is complete to $R_{2,2}$.  If $R_{1,1}$ contains two
adjacent vertices $u,v$, then $\{u, v, x, a_1, a_2\}$ induces an
$F_1$.  So $R_{1,1}$ is a stable set, and similarly $R_{2,2}$ is a
stable set.  Thus $\overline{G}$ is in class $\cal C$.  

Now suppose that $\ell=2$.  Let $z$ be any vertex in
$R\setminus\{x,y\}$.  Suppose that $z$ is anticomplete to $\{a_1,
a_2\}$.  By (\ref{R1}), $z$ is complete to $B$ and has a neighbor $a$
in $A\setminus\{a_1, a_2\}$.  As observed earlier, $a$ is anticomplete
to $\{x,y\}$.  Then $z$ is adjacent to $x$, for otherwise $\{x, z,
a_1, b_1, a\}$ induces a $5$-hole; and similarly $z$ is adjacent to
$b$.  But then $\{x, y, z, a_1, a_2, a\}$ induces a $6$-antihole.
Therefore $z$ has exactly one neighbor in $\{a_1, a_2\}$.  Up to
symmetry, assume that $z$ is adjacent to $a_1$ and not to $a_2$.  If
$z$ is adjacent to $b_1$, then it is also adjacent to $y$, for
otherwise $\{z, a_1, a_2, y, b_1\}$ induces a $5$-hole, and to $x$,
for otherwise $\{z, a_1, x, b_1, y\}$ induces an $F_1$; but then $\{x,
y, a_1, a_2, z\}$ induces an $F_1$.  So $z$ is not adjacent to $b_1$,
and so $z\in R_{2,1}$.  Then $z$ is adjacent to $y$, for otherwise
either $\{z, a_1, a_2, y, b_1, b_2\}$ or $\{z, a_1, a_2, y, b_2\}$
induces a hole, and $z$ is not adjacent to $x$ for otherwise $\{x, y,
a_1, a_2, z\}$ induces an $F_1$.  Then $b_2$ is adjacent to $x$, for
otherwise $\{x, b_1, b_2, z, a_1\}$ induces a $5$-hole, and to $y$,
for otherwise $\{y,b_1, b_2, z, x\}$ induces an $F_1$.  But then
$\{a_1, z, b_2, x, y\}$ induces an $F_1$.  This means that
$R\setminus\{x,y\}=\emptyset$.  If $p\ge 3$, then, as observed
earlier, $\{x,y\}$ is anticomplete to $A\setminus\{a_1, a_2\}$ and
complete to $B$.  It follows that $G$ is in class $\cal C$.  Now
suppose that $p=2$.  Since $x$ and $y$ are $B$-comparable, we may
assume that $N_B(x)\subseteq N_B(y)$.  If $B$ contains two vertices
$b,b'$ that are not adjacent to $x$, then $\{x, a_1, a_2, b, b'\}$
induces an $\overline{F}_1$.  So $B$ has at most one vertex that is
not adjacent to $x$.  If there is such a vertex, then $G$ is in
class~$\cal A$.  If there is no such vertex, then $G$ is in class
$\cal C$.

Now we may assume that any two vertices in $R$ are $A$-comparable and
$B$-comparable.  Since every vertex of $R$ has a neighbor in $A$, some
vertex of $A$ is complete to $R$.  Likewise, some vertex of $B$ is
complete to $R$.  So we may assume that $a_1$ and $b_1$ are complete
to $R$.  If $R$ is not a clique or a stable set, there are three
vertices $x,y,z$ in $R$ that induce a subgraph with one or two edges,
and $\{a_1, b_1, x, y, z\}$ induces an $F_1$ or $F_2$.  Therefore $R$
is a clique or a stable set.
 
Suppose that $R$ is not a clique.  So it is a stable set of size at
least~$2$.  A vertex $a$ in $A\setminus\{a_1\}$ cannot have two
neighbors $x$ and $y$ in $R$, for otherwise $\{a_1, x, y, b, a\}$
induces an $F_1$.  For $k\in\{0,1\}$, let $A_k=\{u\in
A\setminus\{a_1\} \mid u \mbox{ has $k$ neighbors in } R\}$.  So
$A=\{a_1\}\cup A_0\cup A_1$.  Likewise, let $B_k=\{u\in
B\setminus\{b_1\} \mid u \mbox{ has $k$ neighbors in } R\}$.  So
$B=\{b_1\}\cup B_0\cup B_1$.  Since any two vertices in $R$ are
$A$-comparable, some vertex $x$ in $R$ is complete to $A_1$, and
$R\setminus\{x\}$ is anticomplete to $A\setminus\{a_1\}$.  Likewise,
some vertex $y$ in $R$ is complete to $B_1$, and $R\setminus\{y\}$ is
anticomplete to $B\setminus\{b_1\}$.  Suppose that $x=y$.  Consider
any $z\in R\setminus\{x\}$ (recall that $|R|\ge 2$).  Then $z$ is
anticomplete to $(A\setminus\{a_1\})\cup (B\setminus\{b_1\})$, so, by
(\ref{R1}), we have $p=q=2$.  Then $\overline{G}$ is in class $\cal
A$.  Now suppose that we cannot choose $x$ and $y$ equal.  So both
$A_1$ and $B_1$ are not empty and we may assume that $a_2$ is adjacent
to $x$ and not to $y$ and $b_2$ is adjacent to $y$ and not to $x$.  If
there is a vertex $a_0$ in $A_0$, then $\{a_0, a_2, x, y, b_2\}$
induces an $\overline{F}_1$.  So $A_0=\emptyset$.  Likewise
$B_0=\emptyset$.  Thus $G$ is in class $\cal C$.

Now assume that $R$ is a clique.  Since any two vertices of $R$ are
$A$-comparable and $B$-comparable, there is at most one pair $(i,j)$
such that $R_{i,j}\neq \emptyset$, and since $a_1$ and $b_1$ are
complete to $R$, we may assume that $(i,j)=(2,2)$.  Hence $R\setminus
R_0= R_{2,2}$.  Let $R^*=\{x\in R \mid x \mbox{ is complete to } A\cup
B\}$, $R_A=\{x\in R\setminus R^* \mid x \mbox{ is complete to } A\}$
and $R_B=\{x\in R\setminus R^* \mid x \mbox{ is complete to } B\}$.
So $R=R^*\cup R_A\cup R_B\cup R_{2,2}$, and $A\cup R_A$ and $B\cup
R_B$ are cliques.  Since any two vertices in $R$ are $A$-comparable
and $B$-comparable, the bipartite subgraph of $\overline{G}$ induced
by $A\cup R_A\cup B\cup R_B$ is $2K_2$-free; moreover, in that graph
$a_2$ is complete to $B\cup R_B$ and $b_2$ is complete to $A\cup R_A$.
It follows that $\overline{G}$ is in class $\cal B$ (where the four
stable sets are $A\cup R_A$, $B\cup R_B$, $R_{2,2}$ and $R^*$).  This
complete to proof of the theorem.  $\Box$

\medskip

The second condition of Theorem~\ref{thm:main} implies that deciding
whether a graph on $n$ vertices and $m$ edges is (wheel,
antiwheel)-free can be done by brute force in time $O(n^7)$.  So the
problem is polynomially solvable.  Actually, one can use the third
condition of Theorem~\ref{thm:main} to solve the problem in time
$O(n^2)$, as follows: %
\vspace{-.2cm}
\begin{itemize}
\item 
Testing whether $G$ is a $5$-hole or a $6$-hole can be done in time
$O(n)$.
\item 
Testing whether $G$ is a split graph can be done in time $O(m)$; see
\cite{HS81}.
\item 
For each of the classes ${\cal A}$, ${\cal B}$ and ${\cal C}$, testing
membership in the class can be done in time $O(m)$ directly from the
definition of the class (for class $\cal B$, using
Theorem~\ref{thm:bp5}); we omit the details.
\item 
If the above series of tests fails for $G$, one can run it for
$\overline{G}$ in time $O(n^2)$.
\end{itemize}

% \clearpage

\end{document}